\documentclass[11pt,a4paper,leqno]{article}
\usepackage{amsmath}
\usepackage{amssymb}
\advance\topmargin by -2.2cm
\newcommand{\bbr}{I\!\!R}

\newcommand{\bbf}{I\!\!F}

\newcommand{\cala}{{\cal A}}

\newcommand{\calc}{{\cal C}}

\newcommand{\calf}{{\cal F}}

\newcommand{\calt}{{\cal T}}

\newcommand{\barr}{\begin{array}}
\newcommand{\earr}{\end{array}}
\newcommand{\beqq}{\begin{equation}}
\newcommand{\eeqq}{\end{equation}}
\newcommand{\beao}{\begin{eqnarray*}}
\newcommand{\eeao}{\end{eqnarray*}\noindent}
\newcommand{\beam}{\begin{eqnarray}}
\newcommand{\eeam}{\end{eqnarray}\noindent}

\newcommand{\halmos}{\quad\hfill\mbox{$\Box$}}

\newcommand{\la}{\lambda}

\newcommand{\si}{\sigma}
\newcommand{\al}{\alpha}

\newcommand{\Om}{\Omega}

\newcommand{\vep}{\varepsilon}

\newcommand{\wh}{\widehat}
\newcommand{\wt}{\widetilde}

\newcommand{\nto}{n\to\infty}

\begin{document}

{\LARGE 
An extension of the Yamada-Watanabe condition for pathwise\\ 
uniqueness to stochastic differential equations with jumps}
\vskip0.5cm
{\bf R. H\"opfner}, Johannes Gutenberg Universit\"at Mainz \\
\vskip1.0cm

{\bf Abstract: } We extend the Yamada-Watanabe condition for pathwise uniqueness 
to stochastic differential equations with jumps, in the special case where small jumps are summable.  
{\bf MSC: } 60 J 60, 60 J 75 \\


\vskip1.0cm

Results on pathwise uniqueness of solutions for  stochastic differential equations with jumps, driven by Brownian motion $W$ and Poisson random measure $\mu$  
\beam\label{genSDE}
dX_t &=& b(t,X_{t-})\,dt \;+\; \si(t,X_{t-})\,dW_t \\[2mm]
&&+\; \int_{\{|y|\le c\}}f_2(t,X_{t-},y)\,\wt\mu(dt,dy) 
\; +\; \int_{\{|y|>c\}}f_1(t,X_{t-},y)\,\mu(dt,dy)\nonumber  
\eeam
have been obtained under Lipschitz conditions, see Skorohod ([S 65], Chapter 3.2--3.3), Ikeda and Watanabe ([IW 89], Theorem IV.9.1), Bass ([B 04], Theorem 4.1), Protter ([P 05], Chapter V.3). In absence of jumps, Yamada and Watanabe considered   
\beam\label{contSDE}
dX_t \;=\; b(t,X_t)\,dt \;+\; \si(t,X_t)\,dW_t  
\eeam
with non-Lipschitz diffusion coefficient ([YW 71], see Karatzas and Shreve [KS 91] p.\ 291; see also [Y~78]); the example  $\si(t,x) = \sqrt{x\vee 0}\,$ corresponds to Cox-Ingersoll-Ross type diffusions. Yamada and Watanabe proved  pathwise uniqueness for solutions to (\ref{contSDE}) under the following condition (\ref{V2})+(\ref{V3}) on the diffusion coefficient: 
\beqq\label{V2}
|\si(t,x)-\si(t,x')| \;\le\; h(|x-x'|)  \quad\forall\; x,x',t
\eeqq
where  $h:[0,\infty)\to [0,\infty)$ is continuous and nondecreasing, $h(0)=0$, $h(x)>0$ for $x>0$, and 
\beqq\label{V3}
\int_{(0,\vep)} h^{-2}(v)\, dv \;=\; \infty  \quad\mbox{for every $\vep>0$} \;,   
\eeqq
together with a Lipschitz condition on the drift $b(t,x)$. It is interesting to ask for extensions of this Yamada-Watanabe condition to general SDE (\ref{genSDE}).  This question has already been raised by Bass (see the remarks on p.\ 9 in [B 04], and following (1.2) in [B 02]). Theorem 1.1 in [B 02] (see also [BBC 04], and [Z 02]) proves pathwise uniqueness for solutions of 
\beqq\label{SDEbass}
dX_t \;=\; F(X_{t^-})\, dS^\al_t
\eeqq
driven by a symmetric stable process $S^\al$ of index $1<\al<2$, where $F(\cdot)$ is bounded and satisfies a continuity condition (\ref{V2}) with $h(\cdot)$ such that 
\beqq\label{Vbass}
\int_{(0,\vep)} h^{-\al}(v)\, dv \;=\; \infty   \quad\mbox{for every $\vep>0$} \;.    
\eeqq
The proof of this result relies heavily on particular properties of the stable driving process of index $1<\al<2$; a result for case $0<\al<1$ under very weak conditions ([B 02], Theorem 1.2) 
has been revocated subsequently (remarks following theorem 1.1 in [BBC 04]).\\

We prove another type of extension of the Yamada-Watanabe condition for pathwise uniqueness of solutions of (\ref{genSDE}). Our result --~of limited generality since we assume summability of small jumps of the process $X$~-- combines the original Yamada-Watanabe conditions (\ref{V2})+(\ref{V3}) for the diffusive part with a simple Lipschitz condition concerning the small jumps of $\mu$. Big jumps of $\mu$ are irrelevant in view of pathwise uniqueness. As an example, together with a Cox-Ingersoll-Ross type diffusion coefficient, the jump part can be as in (\ref{SDEbass}) with $F(\cdot)$ Lipschitz and $0<\al<1$. \\

This note is organized as follows: i) we recall the general semimartingale setting (as in Jacod and Shiryaev [JS 87] or M\'etivier [M 82]) needed to deal with solutions of equation (\ref{genSDE}); ii) we state the result (theorem 1); iii) we give the proofs together with some related remarks, and point out at which stage the need for summability of small jumps in theorem 1 did arise. \\

\section{Notations, assumptions, result}

On some stochastic basis $(\Om,\cala,\bbf=(\calf_t)_{t\ge 0},P)$, we consider one-dimensional $\bbf$-Brownian motion $W = (W_t)_{t\ge 0}$ and an $\bbf$-Poisson point process $\mu(ds,dy)$ on $(0,\infty){\times}\bbr$. Thus, according to Ikeda and Watanabe ([IW 89], Theorem II.6.3), $\mu$ and $W$ are independent. The random measure $\mu(ds,dx)$ has deterministic intensity 
$$
\wh\mu(ds,dy) \;=\; ds\,\nu(dx)  \quad\mbox{on $(0,\infty){\times}\bbr$}
$$
for some $\si$-finite measure $\nu$ on $\bbr\setminus\{0\}$ satisfying 
\beqq\label{measure_nu}
\int_{\bbr\setminus\{0\}} (y\wedge 1)^2\,\nu(dy) < \infty \;.   
\eeqq
We write $\wt\mu(ds,dy)$ for the compensated random measure 
$$
\wt\mu(ds,dy) \;:=\;  \mu(ds,dy) - \wh\mu(ds,dy)  \quad\mbox{on $(0,\infty){\times}\bbr$}   
$$
and distinguish between small and big jumps of $\mu$ with the help of some $0<c<\infty$; here and below, 'big jump' refers to jumps of the counting process $\left(\mu((0,t]{\times}\{|y|>c\})\right)_{t\ge 0}$. \\

Throughout, we make the following assumptions i)+ii) on the coefficients in equation (\ref{genSDE}): 
 
i) the functions $b(\cdot,\cdot)$ and $\si(\cdot,\cdot)$ are continuous on $[0,\infty){\times}\bbr$, and Yamada-Watanabe conditions  hold (cf.\ [KS 91], p.\ 291): the diffusion coefficient satisfies (\ref{V2}) and (\ref{V3}) above, whereas the drift 
\beqq\label{V1}
|b(t,x)-b(t,x')| \;\le\; K\,|x-x'|  \quad\quad\forall\; x,x',t
\eeqq
is Lipschitz with some constant $K$; 

ii) the functions $(t,x,y) \to f_i(t,x,y)$ are measurable for $i=1,2$; the function $f_2(\cdot,\cdot,\cdot)$ is such that 
$$
\int_{\{|y|\le c\}} f_2^2(t,x,y)\, \nu(dy)  \;<\;  \infty   \quad\quad\forall\; x,t \;;   
$$
whenever we are interested in summability of small jumps of solutions to (\ref{genSDE}), we strengthen this to  
\beqq\label{smalljumpssummabilitycondition}
\int_{\{|y|\le c\}} [\, f_2^2 \vee |f_2|\, ] (t,x,y)\, \nu(dy)  \;<\;  \infty   \quad\quad\forall\; x,t \;.  
\eeqq

\vskip0.3cm
A solution to equation (\ref{genSDE}) is any process $X=(X_t)_{t\ge 0}$ on $(\Om,\cala,P)$ satisfying  iii)--v) below: 

iii) $\,X$ is $\bbf$--adapted and c\`adl\`ag; 

iv) the following process is locally integrable:  
$$
\int_0^t\left\{ |b(t,X_{t-})| + \si^2(t,X_{t-}) \right\} dt \;+\; 
\int_0^t dt \int_{\{|y|\le c\}}f_2^2(t,X_{t-},y)\,\nu(dy) \;,\quad t\ge 0 \;;
$$
whenever (\ref{smalljumpssummabilitycondition}) is assumed, we strengthen this to local integrability of  
$$
\int_0^t\left\{ |b(t,X_{t-})| + \si^2(t,X_{t-}) \right\} dt \;+\; 
\int_0^t dt \int_{\{|y|\le c\}} [\, f_2^2 \vee |f_2|\, ] (t,X_{t-},y)\,\nu(dy) \;,\quad t\ge 0  \;; 
$$
v) the process $X=(X_t)_{t\ge 0}$ has the representation 
\beao
X_t &=& X_0 \;+\;  \int_0^t b(s,X_{s-})\,ds  \;+\;  \int_0^t \si(s,X_{s-})\,dW_s \\
&&+\; \int_0^t\int_{\{|y|\le c\}} f_2(s,X_{s-},y)\,\wt\mu(ds,dy)  
\;+\; \int_0^t\int_{\{|y|>c\}} f_1(s,X_{s-},y)\,\mu(ds,dy) \;. 
\eeao

\vskip0.5cm
These are general conditions needed to deal with solutions of SDE with jumps. In the restricted setting (\ref{smalljumpssummabilitycondition}) where small jumps are summable,  we have the following result.\\

{\bf Theorem 1: } Consider equation (\ref{genSDE})  in case where $f_2$ and $\nu$ satisfy  condition (\ref{smalljumpssummabilitycondition}). 
Together with Yamada-Watanabe conditions (\ref{V2})+(\ref{V3})+(\ref{V1}) on the diffusive part in (\ref{genSDE}) assume 
a Lipschitz condition  
\beqq\label{V4}
\int_{\{|y|\le c\}} |f_2(t,x,y) - f_2(t,x',y)|\,\nu(dy) \;<\; K\, |x-x'|  \quad\forall\; x,x',t  \;.  
\eeqq
Then pathwise uniqueness holds for solutions of equation (\ref{genSDE}). \\

\section{Proofs and some associated results}

We start in the general setting i)--v), without assuming summability of small jumps. The first lemma --~essentially well known as seen from the remarks preceding (3.10) in [B 04], or from p.\ 58 in [S 65]~-- says that big jumps are irrelevant in view of pathwise uniqueness. \\

{\bf Lemma 1: } Let $\calt$ denote the class of $\bbf$-stopping times which are $P$-a.s.\ finite. For every $S\in\calt$, consider the filtration  $\bbf^S:=(\calf_{S+s})_{s\ge 0}$, the $\bbf^S$--Brownian motion $W^S:=(W_{S+s}-W_S)_{s\ge 0}$, and the $\bbf^S$--Poisson point process $\mu^S(ds,dy)$ with intensity $ds\,\nu(dy)$ on $(0,\infty){\times}\bbr$ defined by  $\mu^S(]0,s]{\times}\cdot):=\mu(]S,S{+}s]{\times}\cdot)$, and $\bbf^S$--adapted solutions $X^S$ to equation 
\beam\label{specSDEts}
dX^S_s &=& b(S{+}s,X^S_{s-})\,ds \;+\; \si(S{+}s,X^S_{s-})\,dW^S_s  \\
&&+\; \int_{\{|y|\le c\}}f_2(S{+}s,X^S_{s-},y)\,\wt\mu^S(ds,dy)  \;,\quad s\ge 0 \;. \nonumber
\eeam
If for arbitrary $S\in\calt$ we can prove pathwise uniqueness for equation (\ref{specSDEts}), then pathwise uniqueness holds for equation (\ref{genSDE}).\\


{\bf Proof: } Up to the time dependence in the functions $b$, $\si$, $f_2$, the proof follows [IW 89], p.\ 245.

1) Fix some sequence of constants $(c_m)_m$ with $c_0{:=}c$,  $\,c_m\downarrow 0$ as $m\to\infty$, and $\nu((c_{m+1},c_m])<\infty$ for all $m\ge 0$. 
By basic properties of Poisson random measure ([IW 89], Ch.\ I.9 and II.3),  $1_{\{|y|>c\}}\mu(ds,dy)$ and $1_{\{c_{m+1}<|y|\le c_m\}}\mu(ds,dy)$, $m\ge 0$, are independent random measures. Let $(T_n, Y_n)_n$ denote the sequence of jump times / jump heigths in 
the compound Poisson process 
$\left( \int_0^t \int_{\{|y|>c\}} y\, \mu(ds,dy) \right)_{ t\ge 0}$, 
and write $(S_{m,j})_{j\ge 1}$ for the sequence of jump times of $\left(\mu((0,t]{\times}\{c_{m+1}<|y|\le c_m\}\right)_{t\ge 0}$, for every $m\ge 0$. The graphs (subsets of $[0,\infty){\times}\Omega$, cf.\ [M~82]) 
$$
[[T_n]] \;,\; n\ge 1 \;,\quad [[S_{m,j}]] \;,\; m\ge 0 \;,\; j\ge 1 
$$
are mutually disjoint up to an evanescent set, and support the jumps of $X$. 

2) Let us consider two solutions $\wt X '$, $\wt X ''$ of equation (\ref{genSDE}) with respect to the same pair $(\mu,W)$, starting at time $0$ in the same initial condition $\wt X '_0 = \wt X ''_0$, and let us prove --~under the assumption of the lemma~--  that a.s.\ the paths of $\wt X '$, $\wt X ''$ coincide up to time $\infty$. 

i) First, on the stochastic interval $[[0,T_1[[$, all jumps of $\mu$ are small jumps. As a consequence, before time $T_1$, 
solutions to equation (\ref{genSDE}) are solutions to equation (\ref{specSDEts}) with $S=0$. Hence pathwise uniqueness for equation (\ref{specSDEts}) with $S=0$ yields 
$$
\wt X '\;=\; \wt X '   \quad\mbox{on $\,[[0,T_1[[\,$, a.s.}  \;. 
$$
Since $[[T_1]]$ has (up to an evanescent set) no intersection with $\bigcup_{m,j} [[S_{m,j}]]$ supporting the small jumps of $\mu$, equation (\ref{genSDE}) and step 1) give 
$$
\wt X '_{T_1} \;=\; \wt X '_{T_1^-} + f_1(T_1,\wt X '_{T_1^-},Y_1)  \quad,\quad 
\wt X ''_{T_1} \;=\;  \wt X ''_{T_1^-} + f_1(T_1,\wt X ''_{T_1^-},Y_1)   \;. 
$$
This implies 
\beqq\label{zeit_1}
\wt X '_{T_1} \;=\; \wt X ''_{T_1}  \quad\mbox{a.s.}  
\eeqq
and gives pathwise uniqueness for solutions to (\ref{genSDE}) on the stochastic interval $[[0,T_1]]$. 

ii) Next, consider the solutions $\wt X '$, $\wt X ''$ on the interval $[[T_1,T_2]]$. For $s\ge 0$ and in restriction to the event $\{T_1+s<T_2\}$, representation v) of a solution $\wt X'$ to equation (\ref{genSDE}) gives  
\beam\label{vonzeit1zuzeit2}
\wt X'_{T_1+s}   &=& \wt X'_{T_1} \;+\;  \int_{T_1}^{T_1+s} b(s,\wt X'_{s-})\,ds  \;+\;  \int_{T_1}^{T_1+s} \si(s,\wt X'_{s-})\,dW_s \\
&&+\; \int_{T_1}^{T_1+s}\int_{\{|y|\le c\}} f_2(s,\wt X'_{s-},y)\,\wt\mu(ds,dy) 
\quad\quad\mbox{on}\quad   \{T_1+s<T_2\}  \nonumber
\eeam
since all jumps of $\mu$  on $]]T_1,T_2[[$ are small jumps. The same holds for $\wt X''$ in place of $\wt X'$. Now we put $S:=T_1$ and consider the filtration  $\check{\bbf} :=\bbf^{T_1}$, the $\check{\bbf}$--Brownian motion $\check W:=W^{T_1}$, and the $\check{\bbf}$--Poisson random measure $\check\mu := \mu^{T_1}$, in the notation as above: $\,\check W$ and $\check\mu$ are necessarily independent. 
For $s\ge 0$, put  $\check X '_s:= \wt X '_{T_1+s}$ and $\check X ''_s:= \wt X''_{T_1+s}$. 
Then (\ref{vonzeit1zuzeit2}) shows that before time $\check T_1:=T_2-T_1$  of the first big jump of $\check\mu$, $\;\check X '$ and $\check X ''$ are $\check{\bbf}$-adapted solutions to equation (\ref{specSDEts}) with $S=T_1$, starting from  initial values $\check X '_0$ and $\check X ''_0$ which coincide a.s.\ by (\ref{zeit_1}). By our assumption, pathwise uniqueness holds for equation (\ref{specSDEts}) with $S=T_1$. This show that we have 
$$
\check X '\;=\; \check X ''   \quad\mbox{on $\,[[0,\check T_1[[\,$, a.s.}  \;. 
$$
Changing time back and putting this together with step i), we have  pathwise uniqueness of solutions to (\ref{genSDE}) before time $T_2$. At time $T_2$, we have 
$$
\wt X '_{T_2} \;=\; \wt X '_{T_2^-} + f_1(T_2,\wt X '_{T_2^-},Y_2)  
\;=\; \wt X ''_{T_2^-} + f_1(T_2,\wt X ''_{T_2^-},Y_2) \;=\; 
\wt X ''_{T_2}  \quad\mbox{a.s.}
$$
as above. This gives pathwise uniqueness for solutions to (\ref{genSDE}) on the stochastic interval $[[0,T_2]]$. 

iii) The same argument as in ii) works successively on all intervals $[[T_n,T_{n+1}[[ \,\bigcup\, [[T_{n+1}]]$, $n\ge 1$. Since $T_n\uparrow\infty$, this concludes the proof. \halmos\\


Now we can prove the main result. \\

{\bf Proof of theorem 1: } We have to prove pathwise uniqueness for all equations (\ref{specSDEts}) 
\beao
dX^S_s &=& b(S{+}s,X^S_{s-})\,ds \;+\; \si(S{+}s,X^S_{s-})\,dW^S_s  \\
&&+\; \int_{\{|y|\le c\}}f_2(S{+}s,X^S_{s-},y)\,\wt\mu^S(ds,dy)  \;,\quad s\ge 0 \;.  
\eeao
where $S\in\calt$, according to lemma 1. In equations (\ref{specSDEts}), big jumps are absent.  

I) First, for ease of notation, we consider the particular case  $S=0$ in equation (\ref{specSDEts}). 

1) As in [YW 71] or [KS 91], assumption (\ref{V3}) gives a sequence $a_n \downarrow 0$ such that
$$
a_0=1 \;,\quad \int_{a_n}^{a_{n-1}} h^{-2}(v)\, dv \;=\; n \quad\mbox{for every $n=1,2,\ldots$} \;,   
$$
continuous probability densities $\rho_n(\cdot)$ having support in $(a_n,a_{n-1})$ such that 
$$
\int_{a_n}^{a_{n-1}} \rho_n(v)\, dv \;=\; 1 \quad\mbox{and}\quad   0\le \rho_n(v) \le \frac{2}{n\,h^2(v)} \quad\mbox{for every $n=1,2,\ldots$}  
$$
and $\calc^2$-functions $\psi_n(\cdot)$ on $\bbr$ 
$$
\psi_n(y) \;:=\; \int_0^y \int_0^r \rho_n(v)\,dv\,dr   
\quad\mbox{if $y\ge 0$, and}\quad  \psi_n(y):=\psi_n(-y) \quad\mbox{if $y<0$} \;.  
$$
Then we have  $\psi_n(v)\uparrow |v|$ as $\nto$ for all $v\in\bbr$, and $|\psi_n'(v)|\le 1$ for all $n\ge 1$. 

2) We start without assuming summability of small jumps. By localization, it is sufficient to prove pathwise uniqueness on intervals $[0,T]$ ($T$ deterministic) for solutions $X$ to (\ref{genSDE}) satisfying  
$$
E\left( \int_0^T \left\{ |b(t,X_{t-})|+\si^2(t,X_{t-}) \right\}dt \;+\;  
\int_0^T \int_{\{ |y|\le c \}} f_2^2(t,X_{t-},y)\, \nu(dy)\,dt \right)   \;<\;    \infty  
$$ 
and thus $E(|X_t-X_0|)<\infty$ for all $t\in [0,T]$. Consider two solutions $X^{(1)}$, $X^{(2)}$ to equation  
\beam\label{specSDE}
dX_s \;=\; b(s,X_{s-})\,ds \;+\; \si(s,X_{s-})\,dW_s   
\;+\; \int_{\{|y|\le c\}}f_2(s,X_{s-},y)\,\wt\mu(ds,dy)  
\eeam
with respect to the same pair $(W,\mu)$ and the same initial condition. Then 
$$
D \;:=\; X^{(1)} - X^{(2)}   
$$
has initial value $D_0=0$ and a representation 
\beao
D_t &=& \int_0^t [ b(s,X^{(1)}_{s-}) - b(s,X^{(2)}_{s-}) ]\, ds  \\
&&+\quad  \int_0^t [ \si(s,X^{(1)}_{s-}) - \si(s,X^{(2)}_{s-}) ]\, dW_s \\
&&+\quad  \int_0^t \int_{\{|y|\le c\}} [ f_2(s,X^{(1)}_{s-},y) - f_2(s,X^{(2)}_{s-},y) ]\, \wt\mu(ds,dy) 
\eeao
for $t\in[0,T]$. By Ito formula ([IW 89], p.\ 66) 
\beao
\psi_n(D_t)  &=&   \int_0^t \psi'_n(D_{s-})\, [ b(s,X^{(1)}_{s-}) - b(s,X^{(2)}_{s-}) ]\, ds  \\
&&+\quad  \frac12 \int_0^t \psi_n''(D_{s-})\, [ \si(s,X^{(1)}_{s-}) - \si(s,X^{(2)}_{s-}) ]^2\, ds  \\
&&+\quad  \int_0^t \psi_n'(D_{s-})\, [ \si(s,X^{(1)}_{s-}) - \si(s,X^{(2)}_{s-}) ]\, dW_s \\
&&+\quad \int_0^t \int_{\{|y|\le c\}} [\, \psi_n( D_{s-} + \{f_2(s,X^{(1)}_{s-},y) - f_2(s,X^{(2)}_{s-},y)\} ) - \psi_n( D_{s-} ) \,]\, \wt\mu(ds,dy) \\
&&+\quad \int_0^t \int_{\{|y|\le c\}} \,[\, \psi_n( D_{s-} + \{f_2(s,X^{(1)}_{s-},y) - f_2(s,X^{(2)}_{s-},y)\} ) - \psi_n( D_{s-} ) \\
&&\quad\quad\quad-\;  \{f_2(s,X^{(1)}_{s-},y) - f_2(s,X^{(2)}_{s-},y)\}\,  \psi_n'( D_{s-} ) \,]\, \wh\mu(ds,dy) \;. 
\eeao 
The third and fourth terms on the right hand side are martingales (recall $|\psi_n '(\cdot)|\le1$); all terms on the right hand side are integrable. 

3) The second term on the right hand side of the Ito formula can be treated without any changes as [KS 91], using assumptions (\ref{V2})+(\ref{V3}): we have 
$$
| \si(s,X^{(1)}_{s-}) - \si(s,X^{(2)}_{s-}) |^2  \;\le\;  h^2(|X^{(1)}_{s-}-X^{(2)}_{s-}|)  \;=\; h^2(|D_{s-}|) 
$$
where $\psi_n''=\rho_n$ and $\rho_n\le\frac{2}{n\, h^2}$: for this term,  we have the bound 
\beqq\label{bound5}
\frac12 \int_0^t \psi_n''(D_{s-})\, [ \si(s,X^{(1)}_{s-}) - \si(s,X^{(2)}_{s-}) ]^2\, ds   \;\;\le\;\;  \frac{t}{n}  \;. 
\eeqq
 
\vskip0.2cm
4) Taylor formula with remainder terms written in form 
$$
g(v) \;=\; g(v_0) \;+\;  \sum_{j=1}^{m-1} \frac{g^{(j)}(v_0)}{j\, !}\,(v-v_0)^j \;+\; \int_{v_0}^v \frac{g^{(m)}(r)}{(m-1)\, !}\,(v-r)^{m-1}\, dr 
$$
and short notation 
$$
\zeta(s,y) \;:=\; \{f_2(s,X^{(1)}_{s-},y) - f_2(s,X^{(2)}_{s-},y)\} 
$$
will be used to consider the fifth term 
\beqq\label{T6}
\int_0^t \int_{\{|y|\le c\}} \,[\, \psi_n(\, D_{s-} + \zeta(s,y) \,) - \psi_n( D_{s-} )  \;-\;  \zeta(s,y)\,  \psi_n'( D_{s-} ) \,]\, \wh\mu(ds,dy)  
\eeqq
on the right hand side of the Ito formula. 

i) A first idea is to approximate (\ref{T6}) by 
$$
\int_0^t ds\, \int_{\{|y|\le c\}} \nu(dy)\,  \frac{\psi_n''( D_{s-} )}{2}\, \zeta(s,y)^2 
$$
which equals  
\beam\label{T6-1}
\frac12 \int_0^t ds\, \rho_n(D_{s-})  \int_{\{|y|\le c\}} \nu(dy)\, \{f_2(s,X^{(1)}_{s-},y) - f_2(s,X^{(2)}_{s-},y)\}^2  
\eeam
and would  allow to use  -- instead of our Lipschitz assumption (\ref{V4}) -- a much weaker assumption 
\beqq\label{V4'}
\int_{\{|y|\le c\}} \nu(dy)\, \{f_2(s,x,y) - f_2(s,x',y)\}^2 \;\le\; h^2(|x-x'|) \quad\forall\; x,x',s \;: 
\eeqq
under (\ref{V4'}), the term (\ref{T6-1}) is bounded by $\frac{t}{n}$ in analogy to (\ref{bound5}) above. With this approach however I~was unable to control  remainder terms which involve the heavily fluctuating derivative  $\rho'_n(\cdot)$. 

ii) In the more restrictive setting of summability (\ref{smalljumpssummabilitycondition}) of small jumps as assumed in the theorem, together with the Lipschitz condition (\ref{V4}) on small jumps, remainder terms do not present any difficulty. The localization step in the beginning of 2) now takes the form
$$
E\left( \int_0^T \left\{ |b(t,X_{t-})|+\si^2(t,X_{t-}) \right\}dt \;+\;  
\int_0^T \int_{\{ |y|\le c \}} [\, f_2^2 \vee |f_2|\, ](t,X_{t-},y)\, \nu(dy)\,dt \right)   \;<\;    \infty  
$$
in accordance with (\ref{smalljumpssummabilitycondition}). Write the term (\ref{T6}) as  
\beam\label{T6-2}
&&\int_0^t ds\, \int_{\{|y|\le c\}} \nu(dy)  \left[ \int_{D_{s-}}^{ D_{s-} + \zeta(s,y) } dr\, \psi''_n(r)\, \left(\, D_{s-} + \zeta(s,y) - r \,\right)  \right] \nonumber \\[2mm] 
&&=\quad \int_0^t ds\, \int_{\{|y|\le c\}} \nu(dy)  \left[ \int_0^{\zeta(s,y)} d\tilde r\, \rho_n(D_{s-}+\tilde r)\, (\zeta(s,y)-\tilde r) \right] \;.  
\eeam
For $\la\ge 0$, define a truncated absolute value  $t_\la(\cdot)$ 
by $t_\la(z):=(|z|-\la)\vee 0$ and write the contribution in squared bracketts in (\ref{T6-2}) as 
\beam\label{knackpunkt}
&&1_{\{\zeta(s,y)>0\}} \int_0^\infty d\tilde r\, \rho_n(D_{s-}+\tilde r)\, (\zeta(s,y)-\tilde r)\, 1_{(\tilde r,\infty)}(\zeta(s,y)) \nonumber \\
&&+\quad 1_{\{\zeta(s,y)<0\}} \int_{-\infty}^0 d\tilde r\, \rho_n(D_{s-}+\tilde r)\, (\tilde r-\zeta(s,y))\, 1_{(-\infty,\tilde r)}(\zeta(s,y))  \nonumber \\[2mm]   
&&=\quad \int_{-\infty}^\infty d\tilde r\, \rho_n(D_{s-}+\tilde r)\, t_{|\tilde r|}(\zeta(s,y)) \\[2mm]   
&&\le\quad \int_{-\infty}^\infty d\tilde r\, \rho_n(D_{s-}+\tilde r)\, |\zeta(s,y)| \quad=\quad |\zeta(s,y)|  \nonumber
\eeam 
since $\rho_n(\cdot)$ is a probability density. By definition of $\zeta(s,y)$, we thus obtain the bound  
$$
\int_0^t ds\, \int_{\{|y|\le c\}} \nu(dy)\,  | f_2(s,X^{(1)}_{s-},y) - f_2(s,X^{(2)}_{s-},y) |   
$$
for the fifth term  (\ref{T6}) on the right hand side of the Ito formula, which by (\ref{V4}) is smaller than 
\beqq\label{bound6}
K\, \int_0^t |X^{(1)}_{s-}-X^{(2)}_{s-}|\, ds \;=\;  K\, \int_0^t | D_s |\, ds \;. 
\eeqq

iii) We note the following: as long as $D_{s-}$ might take values in the support $(a_n,a_{n-1})$ of $\rho_n(\cdot)$, we have to account in (\ref{knackpunkt}) above for values of $\tilde r$ which are arbitrarily close to $0$.  

5) Since $|\psi'_n(\cdot)|\le 1$ on $\bbr$, we use assumption (\ref{V1}) to write the first term on the right hand side of the Ito formula as 
\beqq\label{bound1}
\left| \int_0^t \psi'_n(D_{s-})\, [b(s,X^{(1)}_{s-})-b(s,X^{(2)}_{s-})]\, ds \right| \;\le\; 
K\, \int_0^t |X^{(1)}_{s-}-X^{(2)}_{s-}|\, ds \;=\; K\, \int_0^t | D_s |\, ds \;, 
\eeqq
exactly as in [KS 91].

6) Putting together (\ref{bound5})+(\ref{bound6})+(\ref{bound1}) and taking expectations, we deduce from the Ito formula in step~2)
$$
E(\psi_n(D_s)) \;\le\; C_1 \int_0^t E(|D_s|)\, ds \;+\; C_2\frac{t}{n} \;,\quad 0\le t\le T   
$$
for some constants $C_1$, $C_2$, and finish the proof as in [KS 91]: as $\nto$ we have monotone convergence $\psi_n(z)\uparrow|z|$, and thus  
$$
E(|D_s|) \;\le\; C_1 \int_0^t E(|D_s|)\, ds \;,\quad 0\le t\le T  \;;  
$$
the Gronwall lemma gives $E(|D_s|)=0$ for $0\le t\le T$ and concludes part I) of the proof.

II) We prove pathwise uniqueness for equations (\ref{specSDEts}) with arbitrary $\bbf$--stopping times $S\in\calt$.  
Fix $S\in\calt$. If we replace in steps 2)--6) above the functions $b(\cdot,\cdot)$, $\si(\cdot,\cdot)$, $f_2(\cdot,\cdot,\cdot)$ by 
random objects $b(S{+}\cdot,\cdot)$, $\si(S{+}\cdot,\cdot)$, $f_2(S{+}\cdot,\cdot,\cdot)$, 
$\;\bbf$--Brownian motion $W$ by $\bbf^S$--Brownian motion $W^S$ in the notation of lemma 1, 
$\;\bbf$--Poisson point process $\mu$ by the $\bbf^S$--Poisson point process $\mu^S$ as defined in lemma 1, 
and finally $\bbf$--adapted solutions $X^{(i)}$ to (\ref{specSDE}) by $\bbf^S$--adapted solutions $X^{S,(i)}$ to (\ref{specSDEts}), 
then all arguments in steps 2)--6) above will go through exactly as before. The reason is that assumptions (\ref{V2})+(\ref{V3})+(\ref{V1})+(\ref{V4}) allow to vary freely the time argument in the functions  $b(\cdot,\cdot)$, $\si(\cdot,\cdot)$, $f_2(\cdot,\cdot,\cdot)$. This completes the proof of theorem 1. \halmos\\

We add a remark on the case where the heigth of small jumps of the solution process $X$ does not depend on the present state of $X$.   \\

{\bf Proposition 1: } Consider equation (\ref{genSDE}) in case where 
$$
\forall\; t , y\;: \quad  f_2(t,x,y) \;=:\; f_2(t,y) \quad\mbox{does not depend on $x\in\bbr$} \;.  
$$
Then conditions (\ref{V2})+(\ref{V3})+(\ref{V1}) are sufficient for pathwise uniqueness of solutions of equation (\ref{genSDE}). \\

{\bf Proof: } This is a variant of the proof of theorem 1, which does not require the restrictive condition on summability of small jumps used in step 4) of the preceding proof. According to lemma 1, we have to prove pathwise uniqueness for all equations (\ref{specSDEts}) 
\beao
dX^S_s &=& b(S{+}s,X^S_{s-})\,ds \;+\; \si(S{+}s,X^S_{s-})\,dW^S_s  \\
&&+\; \int_{\{|y|\le c\}}f_2(S{+}s,X^S_{s-},y)\,\wt\mu^S(ds,dy)  \;,\quad s\ge 0  
\eeao
where $S\in\calt$. Consider solutions $X^{(1)}, X^{(2)}$ of (\ref{specSDEts}) starting at the same point. In case where the function $f_2(t,x,y)$ does not depend on the space variable $x$, the difference $D := X^{(1)}-X^{(2)}$ is a process with continuous paths     
$$
D_t \;=\; \int_0^t [ b(S{+}s,X^{(1)}_{s-}) - b(S{+}s,X^{(2)}_{s-}) ]\, ds   
\;+\; \int_0^t [ \si(S{+}s,X^{(1)}_{s-}) - \si(S{+}s,X^{(2)}_{s-}) ]\, dW^S_s \;, 
$$
and the Ito formula in step 2) of the preceding proof simplifies to  
\beao
\psi_n(D_t)  &=&   \int_0^t \psi'_n(D_{s-})\, [ b(S{+}s,X^{(1)}_{s-}) - b(S{+}s,X^{(2)}_{s-}) ]\, ds  \\
&&+\quad  \frac12 \int_0^t \psi_n''(D_{s-})\, [ \si(S{+}s,X^{(1)}_{s-}) - \si(S{+}s,X^{(2)}_{s-}) ]^2\, ds  \\
&&+\quad  \int_0^t \psi_n'(D_{s-})\, [ \si(S{+}s,X^{(1)}_{s-}) - \si(S{+}s,X^{(2)}_{s-}) ]\, dW^S_s \;. 
\eeao 
Assuming (\ref{V2})+(\ref{V3})+(\ref{V1}) and localizing as in the beginning of step 2) above,  (\ref{bound5})+(\ref{bound1})  conclude the proof, exactly as in the original Yamada-Watanabe argument for the continuous process (\ref{contSDE}). \halmos \\

\vskip2.0cm
{\bf References: } 


\vskip0.2cm
[B 02]\quad
Bass, R.: Stochastic differential equations driven by symmetric stable processes.\\ 
S\'eminaire de Probabilit\'es (Strasbourg) {\bf 36}, 302--313 (2002). 

\vskip0.2cm
[B 04]\quad
Bass, R.: Stochastic differential equations with jumps. \\
Probability Surveys 1, 1--19 (2004). 

\vskip0.2cm
[BBC 04]\quad
Bass, R., Burdzy, K., Chen, Z.: Stochastic differential equations driven by stable processes for which pathwise uniqueness fails. 
Stoch.\ Proc.\ Appl.\ {\bf 111}, 1--15 (2004). 


\vskip0.2cm
[IW 89]\quad
Ikeda, N., Watanabe, S.: Stochastic differential equations and diffusion processes. \\
2nd ed.\ North-Holland / Kodansha 1989

\vskip0.2cm
[JS 87]\quad
Jacod, J., Shiryaev, A.: Limit theorems for stochastic processes. Springer 1987. 

\vskip0.2cm
[KS 91]\quad
Karatzas, I., Shreve, S.: Continuous martingales and Brownian motion. \\
2nd ed.\ Springer 1991. 

\vskip0.2cm
[M 82]\quad
M\'etivier, M.: Semimartingales. deGruyter 1982. 

\vskip0.2cm
[P 05]\quad
Protter, P.: Stochastic integration and differential equations. 2nd ed.\ Springer 2005

\vskip0.2cm
[S 65]\quad
Skorokhod, A.: Studies in the theory of random processes. Addison-Wesley 1965.  

\vskip0.2cm
[Y 78]\quad
Yamada, T.: Sur une construction des solutions d' \'equations differentielles stochastiques dans le cas non-lipschitzien. 
S\'eminaire de Probabilit\'es (Strasbourg) {\bf 12}, 114--131 (1978). 

\vskip0.2cm
[YW 71]\quad
Yamada, T., Watanabe, S.: On the uniqueness of solutions of stochastic differential equations. 
J.\ Math.\ Kyoto Univ.\ {\bf 11}, 155--167 (1971). 

\vskip0.2cm
[Z 02]\quad
Zanzotto, P.: On stochastic differential equations driven by a Cauchy process and other stable L\'evy motions. 
Ann.\ Prob.\ {\bf 30}, 802-825 (2002).

\vskip2.0cm
\small 
Reinhard H\"opfner

Institut f\"ur Mathematik, Johannes Gutenberg Universit\"at Mainz

D-55099 Mainz, Germany

{\tt hoepfner@mathematik.uni-mainz.de}

{\tt www.mathematik.uni-mainz.de/$\sim$hoepfner}
\vskip0.8cm

\bf{ 2nd revised version 16.09.2009 }


\end{document}